\documentclass[twoside,a4paper,reqno,11pt]{amsart}
\usepackage{amsfonts, amsbsy, amsmath, amsthm, amssymb, latexsym, verbatim, enumerate}
\usepackage{mathrsfs}
\usepackage[top=30mm,right=30mm,bottom=30mm,left=30mm]{geometry}

\usepackage{bm}
\usepackage[pdftex]{color,graphicx}

\makeatletter
\newcommand{\imod}[1]{\allowbreak\mkern4mu({\operator@font mod}\,\,#1)}
\makeatother

\def\hal{\unskip\nobreak\hfill\penalty50\hskip10pt\hbox{}\nobreak

\hfill\vrule height 5pt width 6pt depth 1pt\par\vskip 2mm}

\usepackage{placeins}
\begin{document}
\parskip 1mm

\def\a{\alpha} 
 \def\b{\beta}
 \def\e{\epsilon}
 \def\d{\delta}
  \def\D{\Delta}
\def\r{\rho}
 \def\c{\chi}
 \def\k{\kappa}
 \def\g{\gamma}
 \def\t{\tau}
\def\ti{\tilde}
 \def\N{\mathbb N}
 \def\Q{\mathbb Q}
 \def\Z{\mathbb Z}
 \def\C{\mathbb C}
 \def\F{\mathbb F}
 \def\G{\Gamma}
 \def\go{\rightarrow}
 \def\do{\downarrow}
 \def\ra{\rangle}
 \def\la{\langle}
 \def\fix{{\rm fix}}
 \def\ind{{\rm ind}}
 \def\rfix{{\rm rfix}}
 \def\diam{{\rm diam}}
 \def\M{{\mathcal M}}
\def\QQ{\mathcal{Q}}
\def\veps{\varepsilon}
\def\x{\bar{x}}
\def\y{\bar{y}}
\def\CC{\mathcal{C}}

 \def\rank{{\rm rank}}
 \def\soc{{\rm soc}}
 \def\Cl{{\rm Cl}}

 \def\Sym{{\rm Sym}}
 \def\PSL{{\rm PSL}}
 \def\PSU{{\rm PSU}}
 \def\PGU{{\rm PGU}}
 \def\PGL{{\rm PGL}}
 \def\GL{{\rm GL}}
 \def\SL{{\rm SL}}
 \def\GU{{\rm GU}}
 \def\SU{{\rm SU}}
 \def\Ext{{\rm Ext}}
 \def\E{{\mathcal E}}
 \def\l{\lambda}
 \def\Lie{\rm Lie}
 \def\s{\sigma}
 \def\O{\Omega}
 \def\o{\omega}
 \def\ot{\otimes}
 \def\op{\oplus}
 \def\pf{\noindent {\bf Proof.$\;$ }}
 \def\Proof{{\it Proof. }$\;\;$}
 \def\no{\noindent}
\def\hal{\unskip\nobreak\hfil\penalty50\hskip10pt\hbox{}\nobreak
 \hfill\vrule height 5pt width 6pt depth 1pt\par\vskip 2mm}

 \renewcommand{\thefootnote}{}

\newtheorem{theorem}{Theorem}
 \newtheorem{thm}{Theorem}[section]
 \newtheorem{theor}{Theorem}[subsection]
 \newtheorem{prop}[thm]{Proposition}
 \newtheorem{lem}[theor]{Lemma}
 \newtheorem{lemma}[thm]{Lemma}
 \newtheorem{propn}[theor]{Proposition}
 \newtheorem{cor}[thm]{Corollary}
 \newtheorem{coroll}[theorem]{Corollary}
 \newtheorem{coro}[theor]{Corollary}
 \newtheorem{rem}[thm]{Remark}
 \newtheorem{cla}[thm]{Claim}
 \newtheorem{notation}[thm]{Notation}
 \newtheorem{prob}[thm]{Problem}

\parskip 1mm

\author{Martin R. Bridson}
 \address{M.R. Bridson, Mathematical Institute, University of Oxford, Oxford
OX2 6GG, UK}

\email{bridson@maths.ox.ax.uk}

\author{ David M. Evans}
\address{D.M. Evans, Department of Mathematics,
    Imperial College, London SW7 2BZ, UK}
\email{david.evans@imperial.ac.uk}

\author{Martin W. Liebeck}
\address{M.W. Liebeck, Department of Mathematics,
    Imperial College, London SW7 2BZ, UK}
\email{m.liebeck@imperial.ac.uk}

\author{Dan Segal}
 \address{D. Segal, All Souls College, Oxford OX1 4AL, UK}

\email{dan.segal@all-souls.ox.ac.uk}

\title[Algorithms for finitely presented groups]{Algorithms determining finite simple images of finitely presented groups}

\subjclass[2010]{Primary 20F10; Secondary 20D06 }

    


\begin{abstract} We address the question: for which collections of finite simple groups does there exist an algorithm that determines the images of an arbitrary finitely presented group that lie in the collection? We prove both positive and negative results. For a collection of finite simple groups that contains infinitely many alternating groups, or contains classical groups of unbounded dimensions, we prove that there is no such algorithm. On the other hand, for families of simple groups of  Lie type of bounded rank, we obtain positive results. For example, given any fixed untwisted Lie type $X$ there is an algorithm that determines whether or not an arbitrary finitely presented group has infinitely many simple images isomorphic to   $X(q)$ for some $q$, and if there are finitely many, the algorithm determines them.
\end{abstract} 

\date{\today}
\maketitle

\section{Introduction}\label{intro}

According to a result of the first author and Wilton \cite{BW}, there is no algorithm that, given a finite presentation, determines whether or not the group presented has a nontrivial finite quotient. Consequently there is no algorithm that determines the finite simple images of a finitely presented group. On the other hand, by \cite{J, PF} there {\it is} an algorithm that determines the images of a finitely presented group that are isomorphic to $\PSL_2(q)$ or $\PGL_2(q)$ for some $q$. More precisely, the algorithm determines whether there are finitely or infinitely many values of $q$ such that $\PSL_2(q)$ or $\PGL_2(q)$ is an image, and in the finite case it enumerates all the values of $q$. A similar algorithm  is given in \cite{J2} concerning the images in $\{\PSL_3(q),\,\PGL_3(q),\,\PSU_3(q),\,\PGU_3(q)\}$  of groups with a 2-generator finite presentation. 

Thus the question arises --  for which collections of finite simple groups is there an algorithm that determines the members of the collection that are images of an arbitrary finitely presented group $G$? In this paper we provide answers to this question.

Note that there is an obvious algorithm that determines the finite cyclic images of $G$ (just work with the abelianisation of $G$), so we restrict attention to the non-abelian simple groups.

First, for groups of unbounded rank, we have a negative answer. In the following statement, by the dimension of a classical group we mean the dimension of its natural module.

\begin{theorem}\label{neg} Let ${\mathcal F}$ be a set of finite simple groups that either contains an infinite number of alternating groups, or contains a classical group of dimension $n$ for infinitely many values of $n$.
Then there is no algorithm that, given a finite presentation, does either of the following:
\begin{itemize}
\item[{\rm (i)}] determines whether or not the group presented has at least one quotient in ${\mathcal F}$;
\item[{\rm (ii)}] determines whether or not the group presented has infinitely many quotients in ${\mathcal F}$.
\end{itemize}
\end{theorem}

On the other hand, we are able to prove a positive result for simple groups of fixed rank, reflecting the results on $\PSL_2$, $\PSL_3$  and $\PSU_3$ already mentioned. In order to state our result we need to define some notation. Let $X$ be a fixed untwisted Lie type (e.g. $A_n$ or $E_8$), and for a prime power $q$ let $X(q)$ be the simple group of type $X$ over $\F_q$. If $X(q)$ possesses a graph automorphism of order $d \in \{2,3\}$, denote by $^d\!X(q)$ the corresponding twisted simple group over $\F_q$. (Here we are using the notation $q$ rather than $q^d$, as used for example in \cite[Chapter 14]{car}.) For convenience write $^1\!X(q) = X(q)$. Denote by ${\rm ID}(^d\!X(q))$ the group generated by inner and diagonal automorphisms of $^d\!X(q)$. Recall that ${\rm ID}(^d\!X(q))$ is the fixed point group of a Frobenius endomorphism of the adjoint simple algebraic group of type $X$ (see \cite[Defn. 2.5.10]{GLS}).

\vspace{2mm}
\no {\bf Definition } For a fixed Lie type $^d\!X$ (with $d \in \{1,2,3\}$) define the classes of groups ${\mathcal X}^d$ and ${\mathcal Y}^d$ as follows.
\begin{itemize}
\item[{(i)}] For any $d$, let
\[
{\mathcal Y}^d = \{ Y\;:\; ^d\!X(q) \le Y \le {\rm ID}(^d\!X(q)) \hbox{ for some }q\}.
\]
\item[{(ii)}] For $d=1$, let
\[
{\mathcal X}^1 = \{X(q)\,:\,\hbox{all prime powers }q\}.
\]
\item[{(iii)}] For $d>1$, let ${\mathcal X}^d ={\mathcal Y}^d$.
\end{itemize}

\bigskip 
\no For example, if $X = A_{n-1}\,(n\ge 3)$, there are two $X$-classes:
\[
\begin{array}{l}
{\mathcal X}^1 = \{\PSL_n(q) \,:\, \hbox{all }q \}, \\
{\mathcal X}^2 = \{ Y\;:\; \PSU_n(q) \le Y \le \PGU_n(q) \hbox{ for some }q\}.
\end{array}
\]
Note that apart from types $^d\!X = \,^2\!A_n$,  $^2\!D_n$ and  $^2\!E_6$, all classes ${\mathcal X}^d$ consist of simple groups.

\smallskip 

\begin{theorem}\label{pos} Let $^d\!X$ be a fixed Lie type, and let ${\mathcal X} = {\mathcal X}^d$ be the corresponding class of groups defined above.
 Then there is an algorithm that determines whether or not any given finitely presented group has infinitely many quotients in ${\mathcal X}$; moreover, if there are finitely many such quotients, the algorithm determines them. In particular, the algorithm determines whether or not any given finitely presented group has at least one quotient in ${\mathcal X}$.
\end{theorem}

\no {\bf Remarks }  
\no 1. The proof of Theorem \ref{pos} shows that the same result holds for various subsets of the classes  
${\mathcal X}^d$-- for example, for the collection of groups in ${\mathcal X}^d$ of fixed characteristic. 

\no 2. Theorem \ref{pos} shows, for example, that there is an algorithm that determines whether a finitely presented group has infinitely many quotients $\PSL_n(q)$ ($n$ fixed, $q$ varying); on the other hand, by Theorem \ref{neg} there is no such algorithm for $\PSL_n(q)$ ($q$ fixed, $n$ varying).

\no 3. The algorithms in Theorem \ref{pos} can be taken to be primitive recursive, that is, they do not involve an unbounded search. While the algorithms are in principle explicit, they are not  practical because of our reliance in the proof on results in model theory such as Theorem \ref{T5} and because of the complexity of the formulas involved. 
There is much interest in producing practical algorithms for finitely presented groups -- for example, the algorithms in \cite{J, J2, PF}  have found many applications. Thus -- given a class ${\mathcal X}^d$ for some fixed Lie type $^d\!X$ -- now that we know that there is an algorithm that determines whether or not any given finitely presented group has infinitely many quotients in ${\mathcal X}^d$, a practical such algorithm  would be of great interest.

\bigskip 
Clearly, if a {\it perfect} group $G$ has a quotient $Y$ such that $^d\!X(q) \le Y \le {\rm ID}(^d\!X(q))$, then $Y = \,^d\!X(q)$. Hence for perfect finitely presented groups we have a complete solution to our question:

\begin{coroll}\label{perf}
For a fixed Lie type $^d\!X$, there is an algorithm that determines whether or not any given perfect finitely presented group has infinitely many quotients of the form $^d\!X(q)$, and if there are finitely many, the algorithm determines them.
\end{coroll}

An important ingredient of the proof of Theorem \ref{pos} is the following result, which enables us to distinguish between images in families of simple groups of twisted and untwisted types. We are very grateful to Alex Lubotzky for suggesting this result together with a sketch of a proof using strong approximation, and to Andrei Rapinchuk for much useful information.

\begin{theorem}\label{twist} Let $^d\!X$ be a fixed Lie type, and let $G$ be a finitely generated group.
If $G$ has infinitely many quotients in the class ${\mathcal Y}^d$, then  $G$ also has infinitely many simple quotients 
$X(q) =\, ^1\!X(q)$.
\end{theorem}

For example, applying the result with $^d\!X = \,^2\!A_{n-1}$ shows that a finitely generated group with infinitely many images $\PSU_n(q)$ also has infinitely many images $\PSL_n(q)$. Applying it with $^d\!X = \,^1\!A_{n-1}$ shows  that a finitely generated group with infinitely many images in the set $\{Y: \PSL_n(q) \le Y \le \PGL_n(q) \hbox{ for some }q\}$ also has infinitely many images $\PSL_n(q)$. 


Our algorithmic results would have been more elegant if it were possible to establish
the latter conclusion with unitary 
groups $\PGU_n(q)$, $\PSU_n(q)$ replacing 
$\PGL_n(q)$, $\PSL_n(q)$ (in Theorem \ref{pos} this would enable us to replace our classes of twisted groups ${\mathcal X}^d$ ($d>1$) by just the simple groups in the class). However this is in fact false, as is shown by Corollary \ref{uni}, which for any $n$ divisible by 4, produces a finitely presented group that has infinitely many images in the set $\{Y:\PSU_n(q) \le Y \le \PGU_n(q) \hbox{ for some }q\}$, but has only finitely many images $\PSU_n(q)$.


\section{Proof of Theorem \ref{neg}}

Our proof of Theorem \ref{neg} is based on the following result, taken from \cite{LSres} and \cite{TW}.

\begin{thm}\label{star} {\rm (\cite[Thm. 2.3]{LSres}, \cite[Thm. 2]{TW})} 
There is a function $f:\N\times \N \rightarrow \N$ with the property that for any two finite non-abelian simple groups $T_1,T_2$ of orders $t_1,t_2$, the following hold.
\begin{itemize}
\item[(i)] Every classical simple group of dimension greater than $f(t_1,t_2)$ is generated by a copy of $T_1$and a copy of $T_2$.
\item[(ii)] Every alternating group $A_n$ with $n>f(t_1,t_2)$ is generated by a copy of $T_1$and a copy of $T_2$.
\end{itemize}
\end{thm}

\begin{cor}\label{obv} Let $T_1,T_2$ be finite non-abelian simple groups of orders $t_1,t_2$, and let $n> f(t_1,t_2)$. Then both the alternating group $A_n$, and any classical group $Cl_n(q)$ of dimension $n$, are images of the free product $T_1*T_2$.
\end{cor}

\no {\bf Deduction of Theorem \ref{neg}}

 Let ${\mathcal F}$ be a set of finite simple groups satisfying (i) or (ii) of Theorem \ref{neg}. Suppose for a contradiction that there is an algorithm that determines whether or not any given finitely presented group has at least one quotient in ${\mathcal F}$. 
We shall show that under this assumption there is an algorithm that determines whether or not a finitely presented group has a nontrivial finite image, contrary to the Bridson-Wilton result in \cite{BW}.

Let $G$ be a finitely presented group. By assumption, there is an algorithm that determines whether or not the finitely presented group $G*G$ has an image in ${\mathcal F}$. If the algorithm answers yes, meaning that $G*G$ has an image in ${\mathcal F}$, then $G$ has a nontrivial finite image. If the answer is no, then we claim that $G$ can have no finite non-abelian simple image: for if $G$ had such an image, call it $S$, then $G*G$ would have an image $S*S$, and hence by Corollary \ref{obv}, $G*G$ would have infinitely many images in ${\mathcal F}$, contradicting the output of the algorithm. Finally, we can easily determine whether or not $G$ has a nontrivial finite abelian image by considering the abelianisation of $G$.

Thus we have established that if Theorem \ref{neg} were false, then there would be an algorithm that determines whether or not $G$ has a nontrivial finite image, and as remarked above, this contradicts \cite{BW}. This completes the proof of Theorem \ref{neg}. \hal

\section{Proof of Theorem \ref{pos}}

The proof of Theorem \ref{pos} relies on some representation theory, and also some model theory. We present the required results in Sections \ref{rep} and \ref{model}, and then complete the proof of Theorem \ref{pos} in Sections \ref{sec33}--\ref{detfin}.

\subsection{Actions on adjoint modules}\label{rep}

The main result we shall need is Theorem \ref{adj} below, most of which follows from \cite[Theorem 4]{LSTAMS} and \cite{LSJGT}. 

In the statement of Theorem \ref{adj}, we let $\bar X$ be a simple adjoint algebraic group over an algebraically closed field of characteristic $p>0$, and $F$ a Frobenius endomorphism of $\bar X$ such that $X:= \bar X^F = {\rm ID}(^d\!X(q)) $ is a group of Lie type over $\F_q$, where $q$ is a power of $p$ and $d\in \{1,2,3\}$ (notation of Section \ref{intro}). We say that a subgroup $H$ of $X$ is {\it of the same type} as $X$ if 
$F^*(H) \cong (\bar X^{F_0})'$ for some Frobenius endomorphism $F_0$. Specifically, for such a subgroup $H$, 
one of the following holds (see \cite[Theorem 5.1]{LSroot}):
\begin{itemize}
\item[(i)] $X$ is untwisted, and $^e\!X(q_0) \le H \le {\rm ID}(^e\!X(q_0))$ for some $e \in \{1,2,3\}$ and some subfield $\F_{q_0}$ of $\F_q$;
\item[(ii)] $X= \,^d\!X(q)$ is twisted (so $d=2$ or 3), and $^d\!X(q_0) \le H \le {\rm ID}(^d\!X(q_0))$ for 
some subfield $\F_{q_0}$ of $\F_q$.
\end{itemize}

In each case, we say that $H$ is \textit{over} the field $\F_{q_0}$. 

Recall that $\bar X$ acts on its adjoint module, the Lie algebra $L(\bar X)$. This action is usually irreducible, but there are exceptions in a small number of characteristics; the precise composition factors can be found in \cite[1.10]{LSTAMS}.

\begin{thm}\label{adj} There is an explicitly computable  function $f: \N \rightarrow \N$ such that the following holds. Let $X=  \bar X^F = 
 {\rm ID}(^d\!X(q))$ as above, and let $r$ be the rank of $\bar X$. Suppose $H$ is a subgroup of $X$ such that $H$ acts absolutely irreducibly on every $\bar X$-composition factor of the adjoint module $L(\bar X)$. Then one of the following holds:
\begin{itemize}
\item[{\rm (i)}] $|H| < f(r)$;
\item[{\rm (ii)}] $H$ is of the same type as $X$; 
\item[{\rm (iii)}] $H$ is as in Table $\ref{list}$.
\end{itemize}
\end{thm}

\begin{table}[h]
\caption{Exceptions in Theorem \ref{adj}} \label{list}
\[
\begin{array}{|l|l|}
\hline
X & H \\
\hline
C_n(q),\,p=2 & D_n^\e (q_0),\,D_n^\e (q_0).2\;(\e = \pm, \F_{q_0}\subseteq \F_q) \\
C_4(q),\,p=2 & ^3\!D_4(q_0)\;(\F_{q_0^3}\subseteq \F_q) \\
C_3(q),\,p=2 & G_2(q_0)\;(\F_{q_0}\subseteq \F_q) \\
C_2(q),\,p=2 & H \le Sp_2(q^2).2 \\
^2\!B_2(q),\,p=2 & H \le (q\pm \sqrt{2q}+1).4 \\
A_1(q),\,p=2 & H \le (q\pm 1).2 \\
\hline
\end{array}
\]
\end{table}

\pf Assume that $H$ is not of the same type as $X$. Choose $X_0 \le X$ of the same type as $X$, such that $H\le X_0$ and $X_0$ is minimal with these properties, and write $F^*(X_0) = \,^{d'}\!X(q_0) = (\bar X^{F_0})'$ for some Frobenius endomorphism $F_0$ and $\F_{q_0}\subseteq \F_q$. Let $M$ be a maximal subgroup of $X_0$ such that $H \le M$. Then $M$ is not of the same type as $X$, by the choice of $X_0$.

Suppose $\bar X$ is of exceptional type. Here we apply \cite[Corollary 2]{LSJGT}, which determines the maximal subgroups of $X_0$ that are irreducible on some composition factor of $L(\bar X)$. It follows from this that there is a computable   constant $C$ such that either $|M|<C$, or 
$(\bar X,p)$ is $(F_4,2)$ or $(G_2,3)$; in the latter cases, when $\bar X=F_4$, we have $M=D_4(q_0).S_3$, $C_4(q_0)$ or $B_4(q_0)$, and when $\bar X = G_2$ we have $M = A_2^\e(q_0).2$. However, in none of these cases is $M$ irreducible on both of the composition factors of $L(\bar X)$, contrary to the hypothesis of the theorem. Hence, taking the function $f$ to satisfy $f(r)\ge C$, the result is proved for this case.
\FloatBarrier
Now suppose $\bar X$ is a classical group. We refer to the proof of \cite[Theorem 4.8(i)]{LMT}:  this shows that 
there is a computable function $f_1$ such that either $|M|<f_1(r)$, or there exists a connected $F_0$-stable subgroup $\bar M$ of $\bar X$ such that $M = N_{X_0}(\bar M^{F_0})$, and moreover $M$ stabilizes the Lie algebra $L(\bar M)$. 
Assume $|M|\ge f_1(r)$ (and take the function $f$ to satisfy $f(r)\ge f_1(r)$ for all $r$), so that the latter holds. 
Then $L(\bar X)$ cannot be $\bar X$-irreducible, and so \cite[Prop. 1.10]{LSTAMS} shows that $\bar X$ is as in the following table:
\[
\begin{array}{|l|l|}
\hline
\bar X & \hbox{comp. factor dims. in }L(\bar X) \\
\hline
A_{n-1},\,p|n & n^2-2,\,1 \\
B_n \hbox{ or }C_n,\,p=2 & 2n^2-n-\d,\,2n,\,1^\d\;(\d = {\rm gcd}(n,2)) \\
D_n,\,p=2 & 2n^2-n-\d,\,1^\d\;(\d = {\rm gcd}(n,2)) \\
\hline
\end{array}
\]
If $\bar X = A_{n-1}$ then $\bar M$ must have dimension 1 or $n^2-2$, and it follows easily that $n=2$, $\bar M = T_1$ (a rank 1 torus) and $M = (q_0\pm 1).2$, as in Table \ref{list}. If $\bar X = C_2$ (so that $X_0 = C_2(q_0)$ or $^2\!B_2(q_0)$), then the composition factors of $L(\bar X)$ have dimensions $4,4,1^2$ and we see from the lists of maximal subgroups of these groups in \cite{BHRD} that $H$ is as in Table \ref{list}. And if $\bar X = D_n$ with $n\ge 4$, it is easy to see that $\bar X$ has no connected subgroups of dimension or codimension 1 or 2 satifying the irreducibility hypothesis of the theorem.

So assume now that $\bar X = B_n$ or  $C_n$ with $n\ge 3$ (and $p=2$). The possibilities for the connected subgroups 
$\bar M$ are given by \cite[Theorem 2]{LSINV}: either they belong to one of seven classes of subgroups ${\mathcal C}_i$, or they are simple and act irreducibly on the natural module for $\bar X$. Inspection of the non-simple subgroups in the classes 
${\mathcal C}_i$ shows that none of these has a dimension that is a sum of some of the composition factor dimensions in 
$L(\bar X)$ given above. Hence $\bar M$ is simple. Now \cite[Lemma 7.1]{LSTAMS} shows that the only possibililties for $(\bar X,\bar M)$ are $(C_n,D_n)$, $(B_n,D_n)$ and $(C_3,G_2)$. Hence $X = C_n(q)$ and $M = D_n^\e(q_0).2$ or $G_2(q_0)$ (with $n=3$ in the latter case).  By what we have already proved, the only subgroups of $M$ that satify the irreducibility hypothesis are of the same type as $M$. Hence $H$ is as in Table \ref{list}, and the proof is complete. 
\hal

\no \textbf{Remark } Note that in Theorem \ref{adj}(ii), if $H$ is of the same type as $X$, then $H$ does act absolutely irreducibly on each $\bar{X}$-composition factor of $L(\bar{X})$.

\begin{lemma}\label{spin}
Let $X = C_n(q)$ with $q$ even, $n\ge 3$, and let $V = V(\l_n)$ be a spin module for $X$ over $\F_q$ of dimension $2^n$.
Suppose $H$ is as in the first three rows of Table $\ref{list}$. Then $H$ has a subgroup of index at most $2$ that is not absolutely irreducible on $V$.
\end{lemma}

\pf For the first two rows, observe that $V\downarrow D_n^\e(q)$ is the sum of two half-spin modules for $D_n^\e(q)$. And for the third row, $V \downarrow G_2(q)$ has composition factors of dimensions $6,1,1$. \hal


We also require the following fact about containment between the groups $^d\!X(q)$. 

\begin{thm}\label{inclusion} Let $^d\!X$ be a fixed Lie type, and $p$ a prime, and let $Y$ be a group such that 
$^d\!X(p^e) \le Y \le {\rm ID}(^d\!X(p^e))$, where $^d\!X(p^e)$ is simple. Then $Y$ is a subgroup of 
${\rm ID}(^d\!X(p^f))$ if and only if $e$ divides $f$ and $d$ is coprime to $f/e$. 
\end{thm}

\pf Suppose  $e$ divides $f$ and $d$ is coprime to $f/e$, and write $k = f/e$. 
Let $\bar X$ be the adjoint simple algebraic group over 
$\bar \F_p$ of type $X$, and let $F$ be a Frobenius endomorphism of $\bar X$ such that $\bar X^F = 
{\rm ID}(^d\!X(p^f))$. Then $\bar X^{F^k} = {\rm ID}(^d\!X(p^e))$, and so $Y \le \bar X^{F^k} \le \bar X^F$.

Conversely, suppose $Y \le Z$, where $Y = \,^d\!X(p^e)$, $Z = \,^d\!X(p^f)$. By \cite[Thm. 5.1]{LSroot}, there are Frobenius endomorphisms $F_Y, F_Z$ such that 
\[
Y = O^{p'}(\bar X^{F_Y}),\; Z = O^{p'}(\bar X^{F_Z}).
\]
Now \cite[Lemma 5.3]{LSroot} says that there is a Frobenius automorphism $\tau$ of $Z$ (as defined in \cite[\S5]{LSroot}) such that $Y = O^{p'}(C_Z(\tau))$. The existence of such a Frobenius automorphism  implies that 
 $e$ divides $f$ and $d$ is coprime to $f/e$, as required. \hal

For $d= 1,2,3$ we say that the $d$-part of $e \in \N$ is the largest power of $d$ which divides $e$. 

\begin{cor} \label{cinclusion} Let $^d\!X$ be a Lie type, with corresponding class ${\mathcal Y}^d$. For $1\leq i\leq r$, suppose $K_i \in {\mathcal Y}^d$ is over the field $\F_{p^{e_i}}$, where $e_1,\ldots, e_r$ all have the same $d$-part. Then there is $g \in \N$ such that each $K_i$ is a subgroup of ${\rm ID}(^d\!X(p^g))$.
\end{cor}

\pf Let $g$ be a common multiple of $e_1,\ldots ,e_r$ such that $g/e_i$ is coprime to $d$ for all $i$. Then Theorem \ref{inclusion} implies that $K_i \le {\rm ID}(^d\!X(p^g))$ for all $i$, as required. \hal

\subsection{Model theory}\label{model}

In this section we present some results from model theory that we need.
We work with the usual  language  of rings $L = \{0, 1, +, -, \cdot\}$. `Definable' means `definable without parameters' unless explicitly stated otherwise.

\no {\bf Notation:\/} Denote the set of prime numbers by $\mathcal{P}$ and the set of prime powers by $\mathcal{Q}$. We say that a subset of $\mathcal{Q}$ is \textit{small} if only a finite number of primes divide its elements. A property which holds for all but a small subset of $\mathcal{Q}$ is said to hold for \textit{almost all} elements of $\mathcal{Q}$.

In \cite{Ax}, Ax showed that there is an algorithm to decide, given a closed $L$-formula, whether or not there are infinitely many finite fields in which it holds. The algorithm makes use of G\"odel's Completeness Theorem and involves an unbounded (but terminating) search. A \textit{primitive recursive} algorithm (so not involving unbounded searches) was given by Fried and Sacerdote in \cite{FS}. (Recursive procedures for the other results below had also been given by Ax.) We state their results in the following way:

\begin{thm} \label{T5}  Suppose $\theta$ is a closed $L$-formula. 
\begin{itemize} 
\item[{\rm (i)}] There is a primitive recursive procedure to decide whether or not there is a finite set of primes such that $\theta$ is true in $\F_q$ for all  $q \in \QQ$ which are not divisible by a prime in the finite set. If there is such a finite set, then the procedure also determines it.
\item[{\rm (ii)}] Given a prime $p$, there is a primitive recursive procedure to decide whether or not $\theta$ is true in $\F_{p^e}$ for all but finitely many $e$. If it is, then the procedure also determines the finite set of exceptional exponents.
\item[{\rm (iii)}] There is a primitive recursive procedure to decide whether or not there are infinitely many finite fields $\F_q$ in which $\theta$ is true. If there are only finitely many, then the procedure will determine them.
\end{itemize}
\end{thm}

Note that (i) says that we can determine whether or not $\{q \in \QQ : \F_q \models \theta\}$ is small. As a special case of (i), we can determine whether or not $\theta$ is true in all finite fields.

\smallskip

\pf
(i) This is from Theorem 2 of \cite{FS} (p.228)  with $K = \mathbb{Q}$.

(ii) This is Theorem 3 of \cite{FS} (p.231).

(iii) 
This is stated in \cite{FS} (bottom of page 204), but we give a few more details. First, apply (i) to $\neg\theta$. 
If the procedure finds that the set of all $q\in \QQ$ such that $\F_q \models \neg\theta$ is not small, 
 then there are infinitely many primes $p_1 < p_2 < \ldots$ and exponents $e_i$ such that $\F_{p_i^{e_i}} \models \theta$, and we can stop. Otherwise, (i) identifies a prime $p_0$ such that $\neg\theta$ is true in $\F_{p^e}$ for all primes $p \geq p_0$. For each prime $p < p_0$, we then apply (ii) to $\neg\theta$ to determine whether or not $\theta$ holds in infinitely many $\F_{p^e}$. If it does not, the algorithm determines the finite set of $\F_{p^e}$ in which  $\theta$ holds (cf. Theorem 4 of \cite{FS}).
\hal


We also require the following result, which follows from Theorem 6.7 of \cite{FHJ}.

\begin{thm}\label{316}  Suppose $\theta$ is a closed $L$-formula and $p$ is a fixed prime. Then there are effectively computable constants $N>0$ and $B$ such that, for $e \geq B$, 
\[ \F_{p^e} \models \theta \Leftrightarrow \F_{p^{e+N}} \models \theta.\]
\end{thm}

\medskip

The following is the effective version of the main result of \cite{CDM}. Here, `effective' means computable in a primitive recursive way.

\begin{thm}\label{T7} ({\rm \cite{FHJ}, Theorem 6.4}) For each $L$-formula $\theta(x_1,\ldots, x_n, y_1, \ldots, y_m)$ we can effectively compute a finite set $\{(\theta_i, \mu_i, \veps_i, r_i): i \in I\}$ where:
\begin{itemize}
\item[]$\theta_i(\y)$ is an $L$-formula, $\mu_i$, $\veps_i$ are positive rational numbers,  and $r_i \in \{0,\ldots, n\}$,  
\end{itemize}
such that for each $q \in \QQ$ and $\bar b \in \F_q^m$ there is a unique $i \in I$ with $\F_q \models \theta_i(\bar b)$ and for this $i$, the size $N_q(\bar b)$ of the set $\theta[\F_q, \bar b] = \{ \bar a \in \F_q^n : \F_q \models \theta(\bar a,\bar b)\}$ is either zero or it satisfies
\[ \vert N_q(\bar b) - \mu_iq^{r_i}\vert\leq \veps_iq^{r_i-\frac{1}{2}}.\]
\end{thm}

In the above, we may assume that the formula $(\forall \bar{x})(\neg\theta(\bar{x},\bar{y}))$ is one of the $\theta_i(\bar{y})$. This then picks out those $q, \bar{b}$ where $N_q(\bar{b})$ is zero.


\begin{rem}\label{r35}\rm 
It is well known that, for fixed $d \in \N$,  the field $\F_{q^d}$ and its Frobenius map $x \mapsto x^q$  are first-order definable in the field $\F_q$  in a uniform way, using parameters from $\F_q$. For example, let $\bar{a} = (a_0,\ldots, a_{d-1})$ be a $d$-tuple of parameters representing the coefficients in $\F_q$ of a polynomial $f_{\bar{a}} = a_0 + \ldots + a_{d-1}x^{d-1} + x^d$ which is irreducible over $\F_q$. We identify the field $\F_{q^d}$ with the set $\F_q^d$ with coordinatewise addition and multiplication modulo $f_{\bar{a}}$. The Frobenius map on $\F_{q^d}$ can then be described using a $d^2$-tuple $\bar{b}$ of parameters giving the entries of  a $d \times d$ matrix which gives an automorphism of $\F_{q^d}$ of order $d$. These  conditions on the parameters $\bar{a}, \bar{b}$ can be expressed as $\F_q \models \gamma(\bar{a}, \bar{b})$ where  $\gamma(\bar{x},\bar{y})$ is an appropriate $L$-formula. 

When we make use of this, the actual choice of the parameters $\bar{a}, \bar{b}$ will be irrelevant and we will suppress them in the notation. Moreover, the results for parameter-free formulas in Theorems \ref{T5} and \ref{T7} will still be applicable. For example, the unitary group ${\rm GU}_n(q)$ thought of as a subgroup of ${\rm GL}_n(q^2)$ is definable in $\F_q$ using parameters $\bar{a}, \bar{b}$ as  above (with $d = 2$). A first-order statement $\theta$ about ${\rm GU}_n(q)$ (say, in the language of groups) can be translated into an $L$-formula $\hat{\theta}(\bar{a}, \bar{b})$ (with parameters) whose truth value does not depend on the choice of the particular $\bar{a}, \bar{b}$ which satisfy $\gamma$. Then 
\[
{\rm GU}_n(q) \models \theta \Leftrightarrow \F_q \models \forall \bar{x}\bar{y} (\gamma(\bar{x}, \bar{y}) \to \hat{\theta}(\bar{x}, \bar{y})),
\]
and the decidability results in Theorem \ref{T5} apply to the $L$-formula on the right-hand side here.
\end{rem}

\subsection{Deduction of Theorem \ref{pos}}\label{sec33}

Let $G = \la x_1,\ldots ,x_s \,|\, \r_1,\ldots, \r_k\ra$ be a finitely presented group, with generators $x_i$ and relators $\r_i$. Let $^d\!X$ be a fixed Lie type, and recall the class ${\mathcal Y}^d$ of groups defined in the Introduction. Now let
\[
{\mathcal Y} = \left\{\begin{array}{ll}
\bigcup_e {\mathcal Y^e}, \hbox{ if }d=1 \\
{\mathcal Y}^d, \hbox{ if }d>1
\end{array} \right.
\]
where on the first line, the union is over all $e$  for which the type $^e\!X$ exists. Exclude the case where $^d\!X$ is of Suzuki or Ree type $^2\!B_2$, $^2\!G_2$ or $^2\!F_4$; these will be considered in Section \ref{suzree}.

For a field $K$, let ${\rm ID}(X(K))$ denote the inner-diagonal group of (untwisted) type $X$ over $K$, and let $V = V(K)$ be the adjoint module for ${\rm ID}(X(K))$. For $K = \F_q$, we can regard ${\rm ID}(X(K))$ as a definable subset $H(K)$ of $K^{l^2}$ (where $l = \dim V$) in a uniform way via the action on $V$, as the solutions of polynomial equations defining the connected adjoint algebraic group of type $X$.


If $X = A_n$, $D_n$ or $E_6$ and $K = \F_q$, then $X$ possesses a graph automorphism of order $d \in \{2,3\}$, and the twisted groups ${\rm ID}(^d\!X(q))$ can also be regarded as uniformly definable subsets $H(\F_q)$ of some power of $K$, using Remark \ref{r35}.  This is covered in detail in Section 5.3 of \cite{Ryten}. 

In one exceptional case we vary the above definition of $H(K)$ and $V(K)$. This is the case where $X = C_n\,(n\ge 2)$ and ${\rm char}(K) = 2$. In this case we have ${\rm ID}(X(K)) = X(K)$; let $V_{2n}(K)$ be the natural module, and $S(K)$ a spin module for $X(K)$, and for $K = \F_q$ regard $X(K)$ as a uniformly definable subset $H(K)$ of $K^{l^2}$ via the action on 
$V = V_{2n}(K)\oplus S(K)$, where $l = 2n+2^n$. 

For $H(K)$ defined as above, let $\Phi_0(K)$ be the set of tuples of elements which satisfy the defining relations of $G$. That is,
\[
\Phi_0(K) = \{(g_1,\ldots, g_s) \in H(K)^s : \r_i(g_1,\ldots, g_s) = 1 \mbox{ in $H(K)$ for}\,  i \leq k\}.
\]
This is a definable subset of $K^{sl^2}$. 

For a tuple $\bar g = (g_1,\ldots, g_s) \in \Phi_0(K)$, define $\la \bar g\ra: = \la g_1,\ldots ,g_s\ra$. 
There is an $L$-formula that expresses the condition that the subgroup $\la \bar g\ra$ of $\GL(V)$ is  absolutely irreducible on every $H(K)$-composition factor of $V$. One way of writing down such a formula is as follows. Fix a basis $B$ of $V$. Let $K_1$ be an  extension field of $K$ with $|K_1:K|\le n$ (definable in $K$ as in Remark \ref{r35}),  let $V(K_1)$ be the $K_1$-vector space $V \otimes_K K_1$, and let ${\mathcal W}$ be the collection of subspaces of $V(K_1)$ that are spanned by some proper non-empty subset of $B$. The formula says that for any such $K_1$, there is no pair $h \in GL(V(K_1))$, $W \in {\mathcal W}$  such that each of the conjugates $g_i^h$ stabilizes $W$. 
There is also a formula expressing the condition  $|\la \bar g\ra |>f(r)$, where $f$ is the function in Theorem \ref{adj}.
Hence if we define $\Phi(K)$ to be the set of all tuples $\bar g\in \Phi_0(K)$ such that $\la \bar g\ra$ is absolutely irreducible on every $H(K)$-composition factor of $V$ and has order greater than $f(r)$, then $\Phi(K)$ 
is uniformly definable in $K$. We can regard $\Phi$ here as a parameter-free $L$-formula (with the parameters needed to define the required field extensions absorbed into the formula by making use of Remark \ref{r35}).
In principle, we could write down this formula $\Phi$ explicitly.

In the case where $X = C_n$ and ${\rm char}(K) = 2$, we adjust $\Phi$ to include the condition that subgroups of index at most 4 in $\la \bar g\ra$ arising from subgroups of such index in $G$ are also absolutely irreducible on every composition factor of $V$.

By Theorem \ref{adj} and Lemma \ref{spin}, for any finite field $K = \F_q$, and any tuple $\bar g \in \Phi(K)$, the group $\la \bar g\ra$ is of the same type as $X$, hence is a group in the class ${\mathcal Y}$ over some subfield $\F_{q_0}$ of $\F_q$. Summarising:

\begin{lemma} \label{L3} For each prime power $q$, the finitely presented group $G$ has a homomorphic image in the class ${\mathcal Y}$ over some subfield $\F_{q_0}$ of $\F_q$ if and only if $\F_q \models \exists\x\Phi(\x)$.
\end{lemma}


At this point we apply Theorem \ref{T7} (taking $m=0$) to the parameter-free formula $\Phi(\x)$ and obtain (effectively) a finite set $\{(\theta_i, \mu_i, \veps_i, r_i): i \in I\}$ as in the theorem. Note that as there are no parameter variables in $\Phi$, the $\theta_i$ are $L$-sentences. Informally, this expresses $\Phi(\F_q)$ as a finite union of subsets of dimensions $r_i\,(i\in I)$, for any $q \in \QQ$.

On the other hand, we have the following simple observation.

\begin{lemma} \label{36} For $q \in \QQ$, the group $H(\F_q)$ acts by conjugation on the set $\Phi(\F_q)$. The size of each orbit is  $|H(\F_q)|$.
\end{lemma}

\pf  If $\bar{g} \in \Phi(\F_q)$, then $\langle\bar{g}\rangle \leq H(\F_q)$ is a subgroup of $H(\F_q)$ of the same type as $X$, which has trivial centralizer in $H(\F_q)$ (by Schur's lemma applied to the adjoint representation of $H(\F_q)$). So the stabilizer under conjugation of $\bar{g}$ is trivial.
\hal

It follows that we may assume that the `dimensions' $r_i$ in the above are either $0$ or at least $x$, where $x$ is the dimension of the simple algebraic group of type $X$. We let $I_x = \{i \in I : r_i = x\}$ and $I_y = \{i \in I : r_i > x\}$. Let 
$$\Theta_j = \bigvee_{i \in I_j} \theta_i$$
for $j= x \mbox{ or } y$. 

Define $\QQ (G)$ to be the set of prime powers $q$ such that $G$ has as a homomorphic image a group in the class ${\mathcal Y}$ over $\F_q$. Let $S = \{q : \F_q \models \exists \x \Phi(\x)\}$, and for a prime $p$, let $S_p = \{e : \F_{p^e} \models \exists \x \Phi(\x)\}$. 

\begin{lemma}\label{L9} The set $\QQ(G)$ is infinite if and only if at least one  of the following holds.
\begin{itemize}
\item[(a)] There are infinitely many primes dividing the numbers in the set $S$ (i.e. this set is not small).
\item[(b)] There is a prime $p$ such that $\{e \in \N : \F_{p^e}\models\Theta_y\}$ is infinite.
\item[(c)] $d \geq 2$ and for some $p$, the $d$-part of the elements of $S_p$ is unbounded.
\end{itemize}
\end{lemma}

\pf If the set in (a) is not small, then the number of primes dividing elements of the set $\QQ(G)$ is infinite, by Lemma \ref{L3}. 

Suppose $p$ is as in (b). Then $\vert\Phi(\F_q)\vert$ is (at least) $O(q^{x+1})$ for infinitely many powers $q$ of $p$. As $\vert H(\F_q)\vert$ is $O(q^x)$, it follows that the number of $H(\F_q)$-orbits on $\Phi(\F_q)$ can be arbitrarily large. As $^d\!X(q)$ has boundedly many orbits on subgroups of the same type over a fixed subfield $\F_{q_0}$ (by \cite[Theorem 5.1]{LSroot}), there is no bound on the size of powers of $p$ in $\QQ(G)$.

Suppose $d \geq 2$ and $p$ is as in (c). For $i \in \N$ we can find $e_i \in S_p$ with $d$-part $d^{j_i}$  where $j_1 < j_2 < \ldots$. By Lemma \ref{L3}, there is a subgroup $K_i$ of ${\rm ID}(^d\!X(p^{e_i}))$ which is in the class ${\mathcal Y}$ and which is an image of $G$. This must be  over a subfield $\F_{p^{f_i}}$ of $\F_{p^{e_i}}$, so $^d\!X(p^{f_i}) \leq \, ^d\!X(p^{e_i})$. Thus  by Theorem \ref{inclusion}, $f_i$ divides $e_i$ and $d$ is coprime to $e_i/f_i$. If $i < k$ then 
$e_k$ is divisible by a higher power of $d$ than $e_i$, and so $^d\!X(p^{f_i}) \not\leq \, ^d\!X(p^{e_k})$ by Theorem \ref{inclusion}. Hence $K_k$ is not isomorphic to $K_i$ and it follows that $\QQ(G)$ is infinite.

Conversely suppose that $\QQ(G)$ is infinite. If $\QQ(G)$ is not small then we are in case (a). If $\QQ(G)$ is small then there is a prime $p$ such that $\{ e: p^e \in \QQ(G)\}$ is infinite. 

Suppose first that there is an infinite subset of $\{ e: p^e \in \QQ(G)\}$ whose elements have the same $d$-part. We claim that for infinitely  many powers $q$ of $p$  we have $\F_q \models \Theta_y$.  Assume for a contradiction that this is not the case.

By Corollary \ref{cinclusion} we have an infinite increasing sequence $q_1, q_2, q_3, \ldots$ of powers of $p$ with the property that $\Phi(\F_{q_k})$ contains tuples generating a group in ${\mathcal Y}$ over $\F_q$ for at least $k$ different values of $q$. In particular, the number of orbits of $^d\!X(q_k)$ on this is at least $k$, so 
\[
n_k := \vert \Phi(\F_{q_k})\vert  \geq  k\vert ^d\!X(q_k)\vert.
\]
 But we can also assume that these $\F_{q_k}$ all satisfy the same $\theta_i$ with $i \in I_x$. The estimate in Theorem \ref{T7} gives constants $\alpha, \beta$ such that 
 \[ 
\vert n_k - \alpha q_k^x\vert \leq \beta q_k^{x-\frac{1}{2}}.
\]
 This implies that $n_k / \vert ^d\!X(q_k) \vert$ is bounded above as $k \to \infty$ and so we have a contradiction.
 
 In the remaining case, $\{ e: p^e \in \QQ(G)\}$ has elements of arbitrarily large $d$-part. So $d\geq 2$. If $p^e \in \QQ(G)$ there is $f \in \N$ and a quotient $K$ of $G$ such that $^d\!X(p^f) \geq K \geq \, {\rm ID}(^d\!X(p^e))$. By Theorem \ref{inclusion}, the $d$-parts of $e$ and $f$ are the same. Moreover, $f \in S_p$. Thus (c) holds.  \hal

\begin{cor}\label{primrec} There is a primitive recursive algorithm to decide whether the set $\QQ(G)$ is infinite or not.
\end{cor}

\pf We say how to check the conditions in Lemma \ref{L9}. We 
use the notation and terminology from there. First, by Theorem \ref{T5}(i), we can test whether the set $S$ is small or not using a primitive recursive algorithm. If it is not small, then $\QQ(G)$ is infinite and we can stop. If it is small, then our algorithm gives us an upper bound $p_0$ on the primes which divide its elements. 

For each prime $p \leq p_0$ we decide whether or not (c) of Lemma \ref{L9} holds. We may of course assume that $d \geq 2$. By applying Theorem \ref{316} to the formula $\exists\x \Phi(\x) \wedge \mbox{`characteristic is $p$'}$ we can determine effectively: some $N \in \N$ and a (possibly empty) set of residue classes $a +N\Z$ such that, for sufficiently large $e \in \N$,  we have $e \in S_p$ if and only if $e$ lies in one of these classes. So it remains to decide, given a residue class $a + N\Z$, whether or not it has elements of unbounded $d$-part. But it is easy to show that this happens if and only if the $d$-part of $a$ is greater than or equal to the $d$-part of $N$.

Finally, for  primes $p \leq p_0$ we can check in a primitive recursive way whether there are infinitely many $e$ with  $\F_{p^e}\models\Theta_y$, by Theorem \ref{T5}(ii). \hal

\bigskip
At this point we can complete the proof of Theorem \ref{pos} (apart from the Suzuki and Ree families, postponed until the next subsection). We assume Theorem \ref{twist}, which will be proved in Section \ref{lubseg}. Fix a Lie type $^d\!X$, and let ${\mathcal X}^d$, ${\mathcal Y}^d$ be the corresponding classes of groups defined in the Introduction. 
If $d>1$, then ${\mathcal Y} = {\mathcal Y}^d= 
{\mathcal X}^d$, so Corollary \ref{primrec} shows that there is an algorithm that determines whether or not $G$ has infinitely many quotients in ${\mathcal X}^d$. 
 Now suppose $d=1$, so that ${\mathcal Y} = \bigcup_e {\mathcal Y}^e$.
By Corollary \ref{primrec}, there is an algorithm that determines whether or not $G$ has infinitely many quotients in ${\mathcal Y}$. If there are infinitely many such quotients, then Theorem \ref{twist} shows that there are also infinitely many quotients in ${\mathcal X}^1$. 

This completes the proof of Theorem \ref{pos} apart from the last assertion -- that if the set $\QQ(G)$ is finite, we can compute this set. So suppose $\QQ(G)$ is finite. Note that if $N$ is a positive integer, we can write down a sentence $\Psi_N$ which, in a field $\F_q$, says that there is a tuple in $H(\F_q)$ which satisfies the formula $\Phi(\x)$ and which generates a subgroup of $H(\F_q)$ with at least $N$ elements. As in Lemma \ref{L3}, we have that $\Psi_N$ is true in $\F_q$ if and only if there is a homomorphic image of $G$ in ${\mathcal Y}$ of size at least $N$. So by our assumption, there is some $N$ with the property that $\Psi_N$ is false in all finite fields. Using Theorem \ref{T5}, we can find such an $N$ by testing whether $\Psi_n$ is false in all finite fields for successive values of $n$ until one is found (this involves an unbounded search). This gives a bound $N$ on the size of the groups in $\mathcal{Y}$ which are homomorphic images of $G$, so we can test groups in ${\mathcal Y}$ of size at most $N$  individually (a series of finite computations) to determine whether or not they are homomorphic images of $G$.

In Subsection \ref{detfin}, we give a more complicated argument which avoids the unbounded search in the previous paragraph, hence providing a primitive recursive algorithm for determining $\QQ(G)$ in the case where it is finite.


\subsection{The Suzuki and Ree groups}\label{suzree}

 We are very grateful to Ivan 
 Toma\v si\'c for his help with the material of this section.
The aim is to adapt the above arguments to handle the classes $\mathcal{X}^2$ where $X$ is $B_2$, $G_2$, or $F_4$ (and the characteristic is 2, 3 or 2 respectively). So these are the classes consisting of the Suzuki and Ree groups 
$^2\!B_2(q),\,^2\!G_2(q),\,^2\!F_4(q)$ -- note that for these groups, ${\rm ID}(^2\!X(q)) = \,^2\!X(q)$, so the classes consist only of the simple groups.

The appropriate model-theoretic framework is that of algebraically closed fields with a Frobenius automorphism (\cite{HFrob}), rather than the model theory of finite fields used in the previous subsection.

Let $L^\sigma$ be the first-order language consisting of $L$, the language of rings used previously, together with a unary function symbol $\sigma$. If $q$ is a power of a prime number $p$ and $K$ is an algebraically closed field of characteristic $p$, we denote by $K_q$ the $L^\sigma$-structure $(K; 0,1,+,- , \cdot, \sigma)$  where $\sigma  : K \to K$ is $x \mapsto x^q$. In particular, the fixed field of $\sigma$ is $\F_q$, so the field $\F_q$ is uniformly definable in $K_q$ (that is, using an $L^\sigma$-formula which is independent of $q$). Similarly, for fixed $p$, let $\nu$ denote the Frobenius map $x \mapsto x^p$. If $q = p^n$ then the field $\F_{p^{2n+1}}$ is uniformly interpreted in $K_q$ as the fixed field of $\nu\circ \sigma^2$. Moreover, $(\nu\circ\sigma)^2(x) = \nu(x)$ for all $x \in \F_{p^{2n+1}}$. Thus we have a \textit{uniform} interpretation of the fields $\F_{p^{2n+1}}$ together with a square root of the Frobenius map on these in the structures $K_{p^n}$.

It follows that if $^2\!X(p^{2n+1})$ is a family of Suzuki or Ree groups, then there is a (parameter-free) $L^\sigma$-formula $H$ such that, for an algebraically closed field $K$ of characteristic $p$ (equal to 2 or 3, as appropriate), the solution set $H(K_{p^n})$ is the subgroup $^2\!X(p^{2n+1})$ of the group $X(K)$ (represented by its action on its adjoint module, as in the previous section). 

The algorithm to determine whether our finitely presented group $G$ has infinitely many images in the $X$-class $\mathcal{X}^2$ is then exactly as before, except that we require results for the $L^\sigma$-structures $K_q$ in place of Theorem \ref{T5} (ii) and Theorem \ref{T7}. The following has been communicated to us by Ivan Toma\v si\'c.

The following is Theorem 1.1 of \cite{RT} .

\begin{thm} For each $L^\sigma$-formula $\theta(x_1,\ldots, x_n, y_1, \ldots, y_m)$ there exists a positive constant $C$ and a finite set $D$ of pairs $(d,\mu)$ with $d \in \Z \cup \{\infty\}$ and $\mu \in \Q^+\cup \{\infty\}$ 
such that for each $q \in \QQ$ and $\bar b \in K_q^m$, if  $N_q(\bar b)$ is the size (possibly $\infty$) of the set $\theta[K_q, \bar b] = \{ \bar a \in K_q^n : K_q \models \theta(\bar a,\bar b)\}$, then either $N_q(\bar{b})$ is zero, or there is a pair $(d,\mu) \in D$ with 
 \[ \vert N_q(\bar b) - \mu q^{d}\vert\leq C q^{d-\frac{1}{2}}.
\]

Moreover, for each $(d,\mu) \in D$ there is an $L^\sigma$-formula $\varphi_{d,\mu}(\bar{y})$ with the property that, for each $q$ and $\bar{b} \in K_q^m$, the above estimate for $N_q(\bar{b})$ holds if and only if $K_q \models \varphi_{d,\mu}(\bar{b})$.
\end{thm}


Note that in the above, the obvious conventions about $\infty$ are being used when $N_q(\bar{b})$ is infinite. The paper \cite{RT} makes no claims about determining $C, D$ and the formulas $\varphi_{d,\mu}$. However, Toma\v si\'c informs us that, using the methods of \cite{T_DTGS}, a version of the above can be proved in which these are determined in a primitive recursive way (so providing a replacement for Theorem \ref{T7}). Moreover, in the case where there are no parameter variables $\bar{y}$, the $\varphi_{d,\mu}$ will be such that testing whether or not there are infinitely many $q$ with $K_q \models \phi_{d,\mu}$ can also be done in an effective way (giving a replacement for Theorem \ref{T5}(ii)).

In summary, we have proved the following result, which is Theorem \ref{pos} for the Suzuki and Ree families:

\begin{theorem} Suppose $X = B_2, F_4$ or $G_2$ and the characteristic is $2,2$ or $3$ respectively. Then there is a primitive recursive algorithm to decide whether or not the finitely presented group $G$ has infinitely images in the class $\mathcal{X}^2$.
\end{theorem}

\subsection{Determining finitely many images}\label{detfin}

We continue with the notation of Section \ref{sec33}. In particular $\mathcal{Y}$ is the class defined at the start of that section (or, given the remarks in Section \ref{suzree}, one of the  classes of Suzuki or Ree groups). We will denote ${\rm ID}(^d\!X(q))$ by $H(q)$ (this is consistent with the above notation where  $H(\F_q)$ is the interpretation of ${\rm ID}(^d\!X(q))$ in $\F_q)$. 

We suppose that $\QQ(G)$ is finite and give a primitive recursive algorithm to determine it. From the above, we have a finite list of the possible primes dividing the elements of $\QQ(G)$, so we fix $p$ to be one of these and determine $\QQ_p(G)$, the images in $\mathcal{Y}$ in this characteristic. By Lemma \ref{L3}, we can test whether this is empty or not by applying  Theorem \ref{T5}(ii) to the formula $\neg(\exists\x\Phi(\x))$. 

From now on, we assume that $\QQ_p(G)$ is non-empty.

Let $E = \{ e : \F_{p^e} \models (\exists \x)(\Phi(\x))\}$. Note that this is infinite. By Lemma \ref{L9}, the $d$-part of elements of this set is bounded (and we can compute a bound effectively, as in the proof of Lemma \ref{L9}). Thus it will suffice to show that we can effectively determine the elements of $\QQ_p(G)$ which are defined over fields $\F_{p^e}$ where the $d$-part of $e$ is some fixed power $a$ of $d$. Call this $\QQ_p^a(G)$. In what follows, $q = p^e$ (sometimes with embellishments)  will always denote a power of $p$ with $d$-part $a$. Note that there is a formula $\chi_a$ such that $\F_{p^f} \models \chi_a$ if and only if $f$ has $d$-part $a$.

\begin{notation} \rm  List the groups in $\QQ_p^a(G)$ as $T_1,\ldots, T_m$, where $T_i$ is over the field $\F_{\pi_i}$ for some power $\pi_i$ of $p$. For each $i$ there is $e_i \in \{1,2,3\}$ with $Y_i \leq T_i \leq Y_i^I$, where $Y_i = \, ^{e_i}\!X(\pi_i)$ and $Y_i^I = {\rm ID}(Y_i)$. 
\end{notation}

We shall show how to determine some $q$ such that $H(q)$ contains a copy of each $T_i$. Note that such a $q$ exists, by Corollary \ref{cinclusion}.

Apply Theorem \ref{T7} to the formula $\Phi(\x) \wedge \mbox{`characteristic is $p$'} \wedge \exists\x\Phi(\x)\wedge \chi_a$. 
This gives  a finite set $\{(\theta_i, \mu_i, \veps_i, r_i): i \in I\}$ as in Theorem \ref{T7} such that if $e \in E$ and $q=p^e$, then 
 there is a unique $i \in I$ with $\F_q \models \theta_i$, and for this $i$, 
 \begin{equation} \label{eq2} \vert \vert \Phi(\F_q)\vert  - \mu_iq^{r_i}\vert\leq \veps_iq^{r_i-\frac{1}{2}}. \end{equation}

In the next lemma, $\dim(X)$ denotes the dimension of the simple algebraic group of type $X$.

\begin{lemma}\label{310} Let $i\in I$, and suppose there are infinitely many $e \in E$ with $\F_{p^e} \models \theta_i$.  Then $r_i = \dim(X)$.
\end{lemma}

\pf By assumption, $G$ has finitely many images in $\QQ_p(G)$. Also $H(q)$ has boundedly many orbits on subgroups of the same type over a given subfield (by \cite[Theorem 5.1]{LSroot}). Hence by Lemma \ref{36}, there is a constant $c$ such that 
\[
 \vert H(q)\vert \leq \vert \Phi(\F_{q}) \vert \leq c\vert H(q)\vert.
\]
 As $\vert H(q)\vert = O(q^{\dim(X)})$, the conclusion follows from (\ref{eq2}).
\hal

\begin{rem}\rm \label{311} For $i$ as in Lemma \ref{310}, we may adjust the formula $\theta_i$ so that there exists a rational number $c_i$ with the property that if $(q_j)$ is an increasing sequence of powers of $p$ with $\F_{q_j} \models \theta_i$ for all $j$, then 
\[\lim_{j \to \infty} \frac{\vert H(q_j)\vert}{q_j^{\dim(X)}} = c_i.\]
Assume from now on that this is the case.
\end{rem}

For $j \leq m$, let $\CC_j(q)$ be a system of representatives for the conjugacy classes of subgroups of $H(q)$ which are isomorphic to the group $T_j$ (of course, this could be empty). Let 
\[ t_j(q) = \sum_{T\in \CC_j(q)} \frac{\vert T\vert}{\vert N_{H(q)}(T)\vert}.\]
Let $s_j$ denote the number of orbits (under conjugation) of $T_j$ on the set of $s$-tuples in $T_j$ which satisfy the defining relations of $G$ and which generate $T_j$. 

\begin{lemma} \label{312} We have 
\[ \vert \Phi(\F_{q}) \vert = \vert H(q)\vert \sum_{j=1}^m t_j(q)s_j.\]
\end{lemma}

\pf The number of subgroups of $H(q)$ which are isomorphic to  $T_j$ is equal to $\sum_{T \in \CC_j(q)} \vert H(q)\vert / \vert  N_{H(q)}(T) \vert$. As the number of tuples in $\Phi(\F_{q})$ which generate some $T \cong T_j$ is $\vert T \vert s_j$, we obtain the required result.
\hal

\begin{lemma} There are infinitely many powers $q$ of $p$ with the property that for all $j \leq m$ and for all powers $q'$ of $p$, we have $t_j(q) \geq t_j(q')$.
\end{lemma}

\pf  
Let $\gamma_j = |\CC_j(q)|$. Using the above notation, if $\gamma_j \neq 0$, then $H(q)$ has one conjugacy class of subgroups isomorphic to $Y_j$, each having normalizer isomorphic to $Y^I_j$ in $H(q)$. Thus $\gamma_j$ is equal to the number of conjugacy classes of subgroups isomorphic to $T_j$ in $Y^I_j$. So in this case it follows that 
\[
t_j(q) = \gamma_j\frac{|T_j|}{|N_{Y^I_j}(T_j)|} = \tau_j
\]
depends only on $j$. 

So either $H(q)$ does not contain a copy of $T_j$ and $t_j(q) = 0$, or $t_j(q) = \tau_j$. As there are infinitely many $q$ where $H(q)$ contains a copy of $T_j$ for all $j \leq m$, the result follows. \hal  

With $\tau_j$ as in the proof of the above lemma,  let 
\[ Q = \{ q : t_j(q) = \tau_j \mbox{ for all } j \leq m\}.\]
Let $\alpha = \sum_{j=1}^m \tau_js_j$.

\begin{lemma} \label{314} There is some $i \in I$ with the property that $\F_q \models \theta_i$ for infinitely many $q \in Q$. For such $i$ we have 
\[ \alpha = \frac{\mu_i}{c_i} \geq \frac{\mu_j}{c_j}\]
for all $j \in I$ for which there are infinitely many powers $q'$ of $p$ with $\F_{q'} \models \theta_j$.
\end{lemma}

\pf Suppose $s \in I$ and $(q_u)_{u\in \N}$ is an increasing sequence with $\F_{q_u} \models \theta_s$ and each $t_j(q_u)$ is constant on the sequence. Let $\beta = \sum_j t_j(q_u) s_j$ and $x = \dim(X)$. 
By Lemmas \ref{310} and \ref{312} we have:
\[ \vert \beta - \mu_s \frac{q_u^x}{\vert H(q_u)\vert} \vert < \frac{\veps_s}{q_u^{1/2}}\left(\frac{q_u^x}{\vert H(q_u)\vert}\right).\]
Letting $u \to \infty$, we obtain $\beta = \mu_s/ c_s$.

It follows that if $i,j$ are as in the statement of the Lemma, then $\mu_i/c_i = \alpha$. Moreover,  we may take an infinite sequence of powers $q$ of $p$  where $\F_q$ satisfies $\theta_j$ and where the $t_k(q)$ are constant. Then for such $q$,  
\[ \mu_j/c_j = \sum_k t_k(q)s_k \leq \sum_k \tau_k s_k = \alpha\]
as required. \hal

For notational convenience, take $i = 0$ where $i$ is as in Lemma \ref{314}. Note that by Lemma \ref{314}, this $i$ can be determined explicitly (as the constants $\mu_j$, $c_j$ are effectively computable). We obtain:

\begin{lemma}\label{315} There is an effectively computable constant $A$ such that if $q = p^e$ and $e\geq A$ and $\F_q \models \theta_0$, then $t_j(q) = \tau_j$ for all $j \leq m$. In particular, $H(q)$ contains a copy of each $T_j$.
\end{lemma}

\pf Suppose $\F_q \models \theta_0$. Then, writing $x = \dim(X)$ as usual,
\begin{multline*}\vert \sum_j(\tau_j - t_j(q)) s_j \vert \leq 
c_0^{-1}\left(\left| \mu_0 - \frac{|\Phi(\F_q)|}{q^x}\right| + \left|\frac{|\Phi(\F_q)|}{q^x}- \frac{|\Phi(\F_q)|}{|H(q)|}c_0\right|\right)\\
\leq c_0^{-1}\veps_0q^{-1/2} +c_0^{-1}\frac{|\Phi(\F_q)|}{q^x}\left|1 - \frac{c_0q^x}{|H(q)|}\right|\\
\leq c_0^{-1}\veps_0q^{-1/2} +(\mu_0+1)\left|1 - \frac{c_0q^x}{|H(q)|}\right|,
\end{multline*}
if $q$ is sufficiently large. Note that the left-hand side here is an integer and the right-hand side tends to zero. Thus, as we know explicitly all of the constants involved, we can compute $A$ such that if $q \geq A$, then the left-hand side is zero. As $\tau_j \geq t_j(q)$ for all $j$ and $q$, this implies that $t_j(q) = \tau_j$, as required.
\hal

We can now give the primitive recursive algorithm to determine $\QQ_p^a(G)$, assuming that we have already tested to determine that it is non-empty.
\begin{enumerate}
\item[(1)] Determine the $\{(\theta_i, \mu_i, \veps_i, r_i): i \in I\}$ from  Theorem \ref{T7}, as above.
\item[(2)] Determine for which $\theta_i$ there are infinitely many $e \in E$ with $\F_{p^e} \models \theta_i$ (and discard the rest).
\item[(3)] Adjust the $\theta_i$ as in Remark \ref{311}.
\item[(4)] Take $i$ such that $\mu_i/c_i$ is maximal, and relabel $i=0$.
\item[(5)] Determine the constants $A, B, N$ in \ref{315} and \ref{316} (with $\theta = \theta_0$); let $A_0$ be the maximum of $A, B$.
\item[(6)] Find $q = p^e$ with $A_0 \leq e \leq A_0 + N$ and $\F_{q} \models \theta_0$.
\end{enumerate}
It follows that $H(q)$ contains a copy of each $T_j$: so we can find these by inspection.

\section{Proof of Theorem \ref{twist}}\label{lubseg}

Let $^d\!X$ be a fixed Lie type, and let ${\mathcal Y}^d$ be the corresponding class of groups defined in the Introduction.
Let $G$ be a finitely generated group, and assume that $G$ has infinitely many quotients 
$^d\!\tilde X(q_i) \in {\mathcal Y}^d$ ($i \in \N$), where $^d\!X(q_i) \le \,^d\!\tilde X(q_i)  \le {\rm ID}(^d\!X(q_i))$.
For each $i$, let $\alpha_i: G \to \,^d\!\tilde X(q_i)$ be a surjection. Let $K_i = \bar \F_{q_i}$, and note that 
$^d\!\tilde X(q_i) < X(K_i)$, the adjoint simple algebraic group over $K_i$. Let $V(K_i)$ be the adjoint module for $X(K_i)$, so for each $i \in \N$ we have 
\[
\alpha_i: G \to \,^d\!\tilde X(q_i) < SL(V(K_i)).
\]
Now let ${\mathcal U}$ be a non-principal ultrafilter on $\N$, and take the corresponding ultralimit of the above sequences. This yields 
\[
\alpha: G \to X(K) < SL(V(K)),
\]
where $K = \prod K_i/{\mathcal U}$, an algebraically closed field. We know that $G\alpha_i=\,^d\!\tilde X(q_i)$ acts irreducibly on every $X(K_i)$-composition factor of $V(K_i)$.

Denote the distinguished generators of $G$ by $\bar{x}$ and let $\bar{a}_i$ be the image of these under $\alpha_i$. Let $\bar{a}$ be the image of $(\bar{a}_i)_{i \in \N}$ in the ultra product. Each $\bar{a}_i$ acts irreducibly on the $X(K_i)$-composition factors of $V_i$ and this can be expressed in a first-order way. So by the {\L}os theorem, $\bar{a}$ acts irreducibly on the $X(K)$-composition factors of $V(K)$. The subgroup $A$ generated by $\bar{a}$ is an infinite image of $G$,  and its Zariski closure $\bar A$ is an algebraic subgroup of $X(K)$ which is irreducible on all composition factors of the adjoint module. It follows that $\bar A = X(K)$, apart possibly from the exceptions corresponding to entries in Table \ref{list} of Theorem \ref{adj}, in which $X = C_n$, $p=2$;  these can be excluded by adding the extra first-order condition given in the preamble to Lemma \ref{L3}. Therefore $A$ is Zariski dense in $X(K)$.

Thus we have an algebraically closed field $K$ and a homomorphism $\alpha : G \to X(K)$ with Zariski-dense image. As $G$ is finitely generated, it follows from \cite[Theorem 4.1]{LL} that there is a global field $k_0$ (of the same characteristic as $K$) and a homomorphism $\psi : G \to X(k_0)$ with Zariski-dense image. (Note that in all of this, we are thinking of $X$ as an algebraic group in its adjoint representation.) 

Let $\Gamma = G\psi$. We wish to apply Pink's Strong Approximation Theorem \cite[Theorem 0.2]{Pink}, so we first need to replace the triple $(k_0, X, \Gamma)$ with a \textit{minimal quasi-model} $(k, Y, \Gamma)$, as in \cite[Section 3]{Pink}. So in particular, $k \subseteq k_0$ is a global field, $Y$ is an adjoint algebraic group over $k$ with an isogeny to $X$ over $k_0$ and we can think of $\Gamma$ as a subgroup of $Y(k)$. Then $(k, Y, \Gamma)$ satisfies \cite[Assumption 0.1]{Pink}.  

Denote by $o$ the ring of integers of $k$. 
For a prime $v$ of $k$, denote by $k_v$ the completion of $k$ with respect to $v$, and identify $k$ with a subset of  $k_v$. Let $o_v$ be the ring of integers of $k_v$, and $F_v$ the residue field of $o_v$.

Let $\pi : \hat{Y} \to Y$ denote the simply connected cover of $Y$, and define 
\[
\begin{array}{l}
\Delta = \{ g \in \hat Y(\bar k)\,:\, g\pi \in \Gamma\}, \hbox{ and } \\
\Delta_k = \Delta \cap \hat Y(k).
\end{array}
\]
By Pink's theorem \cite[Theorem 0.2]{Pink}, for almost all primes $v$ of $k$ we have 
\begin{equation} \label{pinksurj}
\Delta_k \cap \hat Y(o_v) \twoheadrightarrow \hat Y(F_v).
\end{equation}

There exists a finite extension $k_1$ of $k$ such that $\Delta \le \hat Y(k_1)$ and $\hat Y$ is split over $k_1$. For a prime $w$ of $k_1$, write ${\mathcal O}_w$ for the ring of integers in the completion $(k_1)_w$. Again we think of $k_1$ as a subset of $(k_1)_w$.

Denote by $T_1$ the set of primes $w$ of $k_1$ such that $\Delta \not \le \hat Y({\mathcal O}_w)$, and define $T = \{ w \cap o : w \in T_1\}$. As $\Delta$ is finitely generated, $T_1$ and $T$ are finite sets. 

Let $V$ be the set of primes of $k$ that split completely in $k_1$. By the Chebotarev density theorem, $V$ has positive density among the primes of $k$. 

Let $v \in V \setminus T$, and let $w$ be a prime of $k_1$ above $v$. Since $v$ splits completely, $k_v \cong (k_1)_w$, so we can identify $k_1$ with a subset of $k_v$. Then $\Delta \le \hat Y(k_v)$, and since $v \not \in T$, it follows that $\Delta \le \hat Y(o_v)$. Hence by (\ref{pinksurj}), for almost all such $v$ we have $\Delta \twoheadrightarrow \hat Y(F_v)$. 
As $\hat Y$ is split over $k_1$, the image $\hat Y(F_v)$ is of untwisted type, and so $\hat Y(F_v)/Z \cong X(q)$, where $q = |F_v|$ and $Z$ is the centre of $\hat Y(F_v)$. Thus $\Delta$ maps onto $X(q)$ for all such $q$. 

We now note that if $S$ is a finite simple image of $\Delta$, then $S$ is also an image of $\Gamma$, and hence of $G$. Indeed, if $N$ is a normal subgroup of $\Delta$ with $\Delta/N = S$, then $ZN/N$ is a central subgroup of $S$ (where $Z = Z(\Delta)$), so $Z \leq N$. But then $S \cong (\Delta/Z)/ (N/Z)$, which is a quotient of $\Gamma$.

This completes the proof of Theorem \ref{twist}. \hal

\section{Almost simple quotients of some unitary groups}\label{unisec}

In this section we exhibit some finitely presented groups that have infinitely many images in the set $\{Y:\PSU_n(q) \le Y \le \PGU_n(q) \hbox{ for some }q\}$, but only finitely many images $\PSU_n(q)$. This explains why, at least with the methods of this paper, in Theorem \ref{pos} we are unable to replace our classes of twisted groups ${\mathcal X}^d$ ($d>1$) by just the simple groups in the class.

Fix $m\geq2$ and an algebraic number field $k$ with $i=\sqrt{-1}\notin k.$
Denote by $^{\ast}$ the non-trivial $k$-automorphism of $k(i)$, and extend
$^{\ast}$ to an automorphism of $\mathbb{C}$. Define
\begin{equation}\label{unidef}
\begin{array}{l}
G(k) =\left\{  g\in\mathrm{GL}_{m}(k(i))\mid g\cdot g^{\ast T}=1\right\}, \\
G_{1}(k)  =G(k)\cap\mathrm{SL}_{m}(k(i)).
\end{array}
\end{equation}
These are the $k$-rational points of algebraic groups $G=\mathrm{U}_{m}
>G_{1}=\mathrm{SU}_{m}$ defined over $k$. 
If $K\subseteq\mathbb{C}$ is a Galois
extension of $k$ with $i\notin K$, then $G_{1}(K)<G(K)$ are defined by
(\ref{unidef}) (with $K$ replacing $k$), while if $i\in K$ we have $\mathrm{GL}%
_{m}(K)\cong G(K)>G_{1}(K)\cong\mathrm{SL}_{m}(K)$, embedded diagonally in
$\mathrm{GL}_{2m}(K).$ 
Suppose $R$ is a subring of $k$ and
$\pi:R\rightarrow F$ is a homomorphism onto a field $F$; then $\pi$ induces a
homomorphism $R[i]\rightarrow F(\sqrt{-1})$ and hence a homomorphism (still
denoted by $\pi$) from $G(R)=G(k)\cap\mathrm{GL}_{m}(R([i])$ into $G(F)$. Here
$G(F)$ is defined by (\ref{unidef}) where $^{\ast}$ is the non-trivial
$F$-automorphism of $F(i)$ if $\sqrt{-1}\notin F$, while if $-1$ is a square in $F$ 
then $G_{1}(F)\cong\mathrm{SL}_{m}(F)$ (see \cite[Section 2.3.3]{PR}).

Let $S$ be a finite set of primes of $k$ and let $R$ be the ring of
$S$-integers of $k$. 
Then $R_{1}\supseteq R[i]\supseteq 2R_{1}$ 
where $R_{1}$ is ring of $S$-integers of $k(i)$.
 For a prime
$\mathfrak{P}$ of $R_{1}$ with $2 \not \in \mathfrak{P}\cap R=\mathfrak{p}$ set
\[
F_{\mathfrak{p}}=R/\mathfrak{p},~q_{\mathfrak{p}}=\left\vert F_{\mathfrak{p}%
}\right\vert ,
\]
\[
E_{\mathfrak{p}}=R_{1}/\mathfrak{P}=(R[i]+\mathfrak{P})/\mathfrak{P}%
\cong\left\{
\begin{array}
[c]{ccc}%
F_{\mathfrak{p}} & \text{if} & q_{\mathfrak{p}}\equiv1~(\operatorname{mod}4)\\
F_{\mathfrak{p}}(i^{\prime}) & \text{if} & q_{\mathfrak{p}}\equiv
3~(\operatorname{mod}4)
\end{array}
\right.
\]
where $i^{\prime}$ is a square root of $-1$.
Let
\[
\begin{array}{l}
P_{1}   =\left\{  \mathfrak{p}\mid q_{\mathfrak{p}}\equiv1~(\operatorname{mod}4)\right\} \\
P_{3}  =\left\{  \mathfrak{p}\mid q_{\mathfrak{p}}\equiv3~(\operatorname{mod}8)\right\} \\
P_{7}   =\left\{  \mathfrak{p}\mid q_{\mathfrak{p}}\equiv7~(\operatorname{mod}8)\right\}  .
\end{array}
\]
Let $\pi_{\mathfrak{p}}:\mathrm{GL}_{m}(R_{1})\rightarrow\mathrm{GL}_{m}(E_{\mathfrak{p}})$ be the quotient map induced by $R_{1}\rightarrow R_{1}/\mathfrak{P}$.
Write $\rho_{q}:\mathrm{U}_{m}(\mathbb{F}_{q})\rightarrow\mathrm{PGU}_{m}(\mathbb{F}_{q})$ for the quotient map modulo the centre.

\begin{prop}\label{uniprop}
Assume\emph{ }that \emph{either} $k$ is not totally real, \emph{or} $S$
contains a prime $\mathfrak{p}$ such that $\left\vert R/\mathfrak{p}\right\vert \equiv1~(\operatorname{mod}4)$. Let
\[
\Gamma=G(R)=G(k)\cap\mathrm{GL}_{m}(R[i]).
\]
Then $\Gamma$ is a finitely presented infinite group and the following
hold.
\begin{itemize}
\item[{\rm (i)}] For almost all primes $\mathfrak{p}$ of $R,$
\[
G_{1}(F_{\mathfrak{p}})\leq\Gamma\pi_{\mathfrak{p}}\leq G(F_{\mathfrak{p}})\leq\mathrm{G\mathrm{L}}_{m}(F_{\mathfrak{p}}).
\]

\item[{\rm (ii)}] If $\mathfrak{p}\in P_{1}$ then $G_{1}(F_{\mathfrak{p}})=
\mathrm{SL}_{m}(F_{\mathfrak{p}}).$

\item[{\rm (iii)}]  If $\mathfrak{p}\in P_{3}\cup P_{7}$ then $G_{1}(F_{\mathfrak{p}})=
\mathrm{SU}_{m}(F_{\mathfrak{p}})\cong\mathrm{SU}_m(q_{\mathfrak{p}})$ and
$G(F_{\mathfrak{p}})=\mathrm{U}_{m}(F_{\mathfrak{p}})\cong\mathrm{U}_m(q_{\mathfrak{p}}).$

\item[{\rm (iv)}]  If $\mathfrak{p}\in P_{3}$ and $4\mid m$ then $\Gamma
\pi_{\mathfrak{p}}\rho_{q_{\mathfrak{p}}}\nleqslant\mathrm{PSU}_{m}(F_{\mathfrak{p}}).$

\item[{\rm (v)}]  Suppose that $m\geq4$. Then with at most finitely many exceptions,
every finite non-abelian simple image of $\Gamma$ is of the form $\Gamma
\pi_{\mathfrak{p}}\rho_{q_{\mathfrak{p}}}$ for some prime $\mathfrak{p}$ of
$R$.

\item[{\rm (vi)}]  $P_{1}$ is infinite and $P_{3}\cup P_{7}$ is infinite. If
$\sqrt{-2}$ $\in k$ then $P_{7}$ is empty.
\end{itemize}
\end{prop}

\begin{cor}\label{uni}
Let $\Gamma=G(R)$ where $R=\mathbb{Z}[\sqrt{-2}]$ and $4\mid m$. Then for
infinitely many primes $p$, $\Gamma$ has an image $X$ satisfying
\[
\mathrm{PSU}_m(p)<X\leq\mathrm{PGU}_m(p),
\]
and all but finitely many of the finite simple images of $\Gamma$ are of the form
$\mathrm{PSL}_m(p^f)$, $f\in\{1,2\}$.
\end{cor}

{\bf Proof of Proposition \ref{uniprop}. }
First observe that $\Gamma$ is an $S$-arithmetic group in the simple algebraic group $G,$ hence is 
finitely presented (\cite[Theorem 5.11]{PR}). That $\Gamma$ is infinite follows
from (i).

(i)  This follows from the Strong Approximation Theorem (\cite[Theorem 7.12]{PR}), 
because $G_{1}=\mathrm{SU}_{m}$ is a connected and simply connected
algebraic group and $G_{1,S}$ contains $G_{1}(\mathbb{C})\cong\mathrm{SL}_{m}(\mathbb{C})$ if 
$k$ is not totally real, or $G_{1}(k_{\mathfrak{p}})\cong\mathrm{SL}_{m}(k_{\mathfrak{p}})$ where $\mathfrak{p} \in S$ splits in $k(i)$.

(ii) If $\mathfrak{p}\in P_{1}$ then $-1$ is a square in
$F_{\mathfrak{p}}$. Hence (ii).

(iii) If $\mathfrak{p}\in P_{3}\cup P_{7}$ then $-1$ is not a square
in $F_{\mathfrak{p}}$. Hence (iii).

(iv) Assume now that $\mathfrak{p}\in P_{3}$ and write
$q=q_{\mathfrak{p}}=\left\vert F_{\mathfrak{p}}\right\vert $, so
$q\equiv3~(\operatorname{mod}8)$ and $\left\vert E_{\mathfrak{p}}\right\vert
=q^{2}$. Thus for $x\in E_{\mathfrak{p}}$ we have $x^{\ast}=x^{q}$.

Let $g=\mathrm{diag}(-1,1,\ldots,1)\in\mathrm{GL}_{m}(k(i))$. Then $g\in
\Gamma$. Suppose that $g\pi_{\mathfrak{p}}\rho_{q_{\mathfrak{p}}}
\in\mathrm{PSU}_{m}(F_{\mathfrak{p}})$. Then there exist $h\in\mathrm{SU}_{m}(F_{\mathfrak{p}})$ and 
$\mu\in E_{\mathfrak{p}}$ such that $\mu \mathbf{1}_{m}\in\mathrm{U}_{m}(F_{\mathfrak{p}})$ and 
$\mu h=\mathrm{diag}(-1,1,\ldots,1)$. Thus
\[
\begin{array}{l}
\mu^{q+1}  =\mu\mu^{\ast}=1, \\
-1   =\det(\mu h)=\mu^{m}.
\end{array}
\]
If $4\mid m$ it follows that $8\mid o(\mu)\mid q+1$ where $o(\mu)$ is the
order of $\mu$ in $E_{\mathfrak{p}}^{\ast}$, a contradiction.

(v) Now assume that $m\geq4.$ Then $\Gamma$ has the (weak) congruence
subgroup property (see \cite{R} or \cite{T}). This implies that if $N\vartriangleleft\Gamma$
and $\Gamma/N$ is finite then $N\geq\lbrack H,\Gamma]$ for some principal
congruence subgroup $H$. If also $\Gamma/N$ is simple and non-abelian it
follows that $N\geq H$.

Now $H$ is the kernel of the map $\pi_{I}$ induced by $R\rightarrow R/I$ for
some ideal $I\neq0$ of $R$. Say $I= {\textstyle\bigcap}\mathfrak{p}_{j}^{e_{j}}$. 
Each non-abelian simple quotient of $\Gamma\pi_{I}$
is then a quotient of $\Gamma\pi_{\mathfrak{p}_{j}}$ for some $j$. So
$\Gamma/N$ is a quotient of $\Gamma\pi_{\mathfrak{p}}$ for some $\mathfrak{p}$.

It follows from (i), (ii) and (iii) that for almost all $\mathfrak{p}$, the group
$\Gamma\pi_{\mathfrak{p}}$ has at most one simple quotient and that quotient
is $\Gamma\pi_{\mathfrak{p}}\rho_{q_{\mathfrak{p}}}$. The claim follows.

(vi) Observe that $\mathfrak{p}\in P_{1}$ if
$\mathfrak{p}$ splits in $k(i)$, and $\mathfrak{p}\in P_{3}\cup P_{7}$ if
$\mathfrak{p}$ is inert in $k(i)$. Hence both of these sets are infinite, by
(an elementary case of) Chebotarev's Theorem (see e.g. \cite[Theorem 1.2]{PR}).

Now assume that $\sqrt{-2}$ $\in k$, $\sqrt{-1}$ $\notin k.$ Then
$R_{0}=\mathbb{Z}[\sqrt{-2}]\subseteq R$. Suppose $\mathfrak{p}\in P_{7}$, and put
$\mathfrak{p}_{0}=\mathfrak{p}\cap R_{0}$ and $q_{0}=\left\vert R_{0}/\mathfrak{p}_{0}\right\vert $. 
Then $-2$ is a square and $-1$ is not a square
in the field $R_{0}/\mathfrak{p}_{0}$. This implies that $q_{0}\equiv
3~(\operatorname{mod}4),$ hence that $q_{0}=p$ is prime, and then
\[
1= {-2 \choose p} = -{2\choose p}
\]
so $p\equiv\pm3~(\operatorname{mod}8)$. But $q_{\mathfrak{p}}$ is a power of
$p$ so $q_{\mathfrak{p}}\equiv1$ or $\pm3~(\operatorname{mod}8)$,
a contradiction. Hence $P_7$ is empty.

This completes the proof. \hal

\end{document}